\documentclass[12pt,reqno]{article}

\usepackage[usenames]{color}
\usepackage{amssymb}
\usepackage{amsmath}
\usepackage{amsthm}
\usepackage{amsfonts}
\usepackage{amscd}
\usepackage{graphicx}
\usepackage{booktabs}
\usepackage[colorlinks=true,
linkcolor=webgreen,
filecolor=webbrown,
citecolor=webgreen]{hyperref}

\definecolor{webgreen}{rgb}{0,.5,0}
\definecolor{webbrown}{rgb}{.6,0,0}

\usepackage{color}
\usepackage{fullpage}
\usepackage{float}

\usepackage{graphics}
\usepackage{latexsym}

\allowdisplaybreaks

\begin{document}

\theoremstyle{plain}
\newtheorem{theorem}{Theorem}
\newtheorem{proposition}[theorem]{Proposition}
\newtheorem{corollary}[theorem]{Corollary}
\newtheorem{lemma}[theorem]{Lemma}
\theoremstyle{definition}
\newtheorem{example}[theorem]{Example}
\newtheorem*{remark}{Remark}

\numberwithin{equation}{section}
\numberwithin{theorem}{section}

\begin{center}
\vskip 1cm
{\large\bf Double sums involving quotients of binomial coefficients}

\vskip 1cm

{\large
Kunle Adegoke \\
Department of Physics and Engineering Physics \\ Obafemi Awolowo University \\ 220005 Ile-Ife, Nigeria \\
\href{mailto:adegoke00@gmail.com}{\tt adegoke00@gmail.com} \\
ORCID: 0000-0002-3933-0459

\vskip .25 in

Robert Frontczak \\
Independent Researcher \\ 72764 Reutlingen,  Germany \\
\href{mailto:robert.frontczak@web.de}{\tt robert.frontczak@web.de}\\
ORCID: 0000-0002-9373-5297

\vskip .25 in

Karol Gryszka \\
Institute of Mathematics \\ University of the National Education Commission, Krakow \\
Podchor\c{a}\.{z}ych 2, 30-084 Krak{\'o}w, Poland \\
\href{mailto:email}{\tt karol.gryszka@uken.krakow.pl} \\
ORCID: 0000-0002-3258-3330

\vskip .25 in
}

\end{center}

\vskip .2 in

\begin{abstract}
This paper is the third part of our detailed investigation of double sums where the inner sum is binomial but incomplete. 
We prove many identities for these double sums involving quotients of binomial coefficients and additional parameters. 
In addition, we evaluate some quadratic incomplete sums using Legendre polynomials.
\end{abstract}

\noindent 2020 {\it Mathematics Subject Classification}:
Primary 05A10, Secondary 11B65, 33B15, 33C05.

\noindent \emph{Keywords:}
Sum, binomial coefficient, harmonic number, Legendre polynomial.

\section{Introduction}

In our paper \cite{Adegoke3} we have initiated a detailed analysis for double sums of the form
$$\sum_{k=0}^n x^k \left (\sum_{j=0}^k \binom{n}{j} y^j\right ),$$
where $x$ and $y$ are two complex parameters and where the inner sum is binomial but incomplete. 
We have evaluated a range of double sums with Fibonacci numbers $F_n$, Bernoulli numbers $B_n$ and Catalan numbers $C_n$ in closed form. 
In our second paper on the topic \cite{Adegoke4} we have extended the analysis to double sums associated with binomial transforms
and derived, among other things, new combinatorial identities for theses sequences and their polynomial generalizations. \\

The present paper is the final part of our exploration in this direction. Our work was highly inspired by the article written 
by Stenlund and Wan \cite{Stenlund} where the authors analyzed classes of double sums involving ratios of binomial coefficients. 
These sums are of the form
\[ \sum_{k=0}^{n} \sum_{j=0}^k \frac{\binom{f_1(n)}{j}}{\binom{f_2(n)}{k}}\,c^{j-k}, \]
where $c$ is a parameter and $f_1, f_2$ are certain functions of $n$. They derived some of their results using the following two identities
\begin{equation*}
\sum_{j=m}^n \frac{\Gamma(j+a)}{\Gamma(j+b+1)} = \frac{\Gamma(m+a)}{(b-a) \Gamma(m+b)} - \frac{\Gamma(n+a+1)}{(b-a)\Gamma(n+b+1)} 
\quad (n \ge m)
\end{equation*}
and
\begin{equation*}
\sum_{j=0}^{k} \binom{n}{j} x^j = \frac{x^{k+1}}{x+1} \binom{n}{k} \,{}_2F_1\left(1,n+1;n+1-k;\tfrac{1}{x+1}\right),
\end{equation*}
where ${}_2F_1$ is the Gaussian hypergeometric function 
\begin{equation*}
{}_2F_1(a,b;c;z) = \sum_{n=0}^\infty \frac{(a)_n (b)_n}{(c)_n} \frac{z^n}{n!},
\end{equation*}
and $(a)_n$ is the Pochhammer symbol given by $(a)_n = \Gamma(a+n)/\Gamma(a)$, and $a,b,c$ are complex numbers. 
Two particular nice results are 
\begin{equation}\label{Sten1}
\sum_{k=0}^n\sum_{j=0}^k \frac{\binom{n+1}{j}}{\binom{n}{k}} = \frac{n+1}{2} \sum_{k=1}^{n+1} \frac{2^k}{k}
\end{equation}
and the alternating identity
\begin{equation}\label{Sten2}
\sum_{k=0}^n\sum_{j=0}^k \frac{\binom{n+1}{j}}{\binom{n}{k}} (-1)^k = \frac{n+1}{2n+4}\left (1+(-1)^n(2^{n+2}-1)\right).
\end{equation}
In some other cases the authors succeeded to express these double sums in terms of the harmonic numbers $H_n$ and odd harmonic numbers $O_n$, respectively (see Corollary 6 in \cite{Stenlund}), where, as usual, these numbers are defined by
\begin{equation*}
H_0 = 0, \quad H_n = \sum_{m=1}^n \frac{1}{m} \qquad \text{and} \qquad O_0 = 0, \quad O_n = \sum_{m=1}^n \frac{1}{2m-1}.
\end{equation*}
For a complex number $z$ that is not a negative integer, we use the generalized harmonic number 
\[ H_z=\psi(z+1)+\gamma, \] 
where $\psi$ is the digamma function and $\gamma$ is Euler's constant. We also express some sums with the Catalan numbers $C_n$, which are defined by 
$$C_n=\frac{1}{n+1}\binom{2n}{n}.$$

In this article we consider similar forms of sums. However, our approach is different and is not based on hypergeometric functions. 
Instead, we make use of the ideas developed by Zhang \cite{Zhang1,Zhang2}. Zhang, when working on generalizations of Calkin's identity, derived the following result 
\cite[Theorem 1.2, Part 1]{Zhang1}: Let $f_k$ and $g_k$ be two real (or complex) number sequences. Then
\begin{equation}\label{Zhang_id}
\sum_{k=0}^n f_k \left (\sum_{j=0}^k \binom{n}{j}\,g_j \right ) = S_f(n+1)\cdot H_g(n) - P(n,0),
\end{equation}
with
\begin{equation*}
S_f(j) = \sum_{k=0}^{j-1} f_k, \quad H_g(n) = \sum_{j=0}^n \binom{n}{j}\,g_j, \quad \text{and} \quad 
P(n,k) = \sum_{j=k}^n \binom{n}{j}\,g_j\,S_f(j).
\end{equation*}
Zhang's theorem has a Part 2 dealing with quadratic sums which we will discuss later. Based on Zhang's approach we will complement the analysis of Stenlund and Wan, showing that, for instance, 
\begin{equation}\label{showcase1}
\sum_{k=0}^n\sum_{j=0}^k  \frac{\binom{n}{j}}{\binom{n}{k}} (-1)^{k+j} = \frac{n+1}{2}, \quad n\geq 1,
\end{equation}
and 
\begin{equation}\label{showcase2}
\sum_{k=0}^n\sum_{j=0}^k \frac{\binom{n}{j}}{\binom{n+1}{k}} (-1)^{k} 
= \frac{1}{n+3}\left ((-1)^n 2^{n} + \frac{14n^2+16n+14+2(-1)^n}{8(n+1)}\right), \quad n\geq 0,
\end{equation}
where in \eqref{showcase2} the ratio of binomial coefficients is the "inverse ratio" in the double sum \eqref{Sten2}. 
Both, \eqref{showcase1} and \eqref{showcase2}, are special instances of the following more involved results containing an additional 
parameter $y$ and presented below as corollaries:
\begin{align}
& \sum_{k=0}^n\sum_{j=0}^k \frac{\binom{n}{j}}{\binom{n}{k}} (-1)^k y^j \nonumber \\
& \quad = \begin{cases}
  \frac{n+1}{2} + (-1)^n (n+1)\frac{0^n}{n+2}, &\text{if $y=-1$}; \\ 
  \frac{1}{n+2}\left ( (-1)^n (n+1) (1+y)^n + \frac{(-1)^n y^{n+2} + n(1+y)+2y+1}{(1+y)^2}\right ),&\text{if $y\neq -1$}; 
  \end{cases}
\end{align}
and
\begin{align}
& \sum_{k=0}^n\sum_{j=0}^k \frac{\binom{n}{j}}{\binom{n+1}{k}} (-1)^k y^j \nonumber \\
& \quad = \begin{cases}
  \frac{n^2+5n+6}{3(n+3)} + (-1)^n \frac{0^n}{n+3}, &\text{if $y=-1$}; \\ 
  \frac{1}{n+3}\left ( (-1)^n (1+y)^n + \frac{n^2 (1+y)^{2} + n(5y^2+8y+3)+6y^2+6y+2-2(-y)^{n+3}}{(n+1)(1+y)^3}\right ),&\text{if $y\neq -1$}, 
  \end{cases}
\end{align}
where $y$ is an arbitrary real or complex number. The quantities $(-1)^n (n+1)\frac{0^n}{n+2}$ and $(-1)^n \frac{0^n}{n+3}$ are correcting 
factors to make the identities also hold for $n=0$. We will also analyze a range of similar sums of alternating and non-alternating nature.
Finally, we evaluate some quadratic incomplete sums using Legendre polynomials.

\section{An alternating double sum}

\begin{lemma}\label{lem1}
For integers $c,d,i$ and $n$ with $0\leq i \leq n$ and $0\leq d \leq c$ we have
\begin{equation}
\sum_{k=0}^i \frac{(-1)^k}{\binom{n+c}{k+d}} = \frac{n+c+1}{n+c+2}\left (\frac{(-1)^i}{\binom{n+c+1}{i+d+1}} + \frac{1}{\binom{n+c+1}{d}} \right ).
\end{equation}
In particular, 
\begin{equation}\label{lem1_id2}
\sum_{k=0}^i \frac{(-1)^k}{\binom{n+c}{k}} = \frac{n+c+1}{n+c+2}\left (\frac{(-1)^i}{\binom{n+c+1}{i+1}} + 1 \right ).
\end{equation}
\end{lemma}
\begin{proof}
We use the Beta integral approach to prove the identity. We have
\begin{align*}
\sum_{k=0}^i \frac{(-1)^k}{\binom{n+c}{k+d}} &= \sum_{k=0}^i (n+c+1) (-1)^k\int_0^1 t^{k+d} (1-t)^{n+c-d-k}\,dt \\
&= (n+c+1) \int_0^1 t^d (1-t)^{n+c-d}\left (\sum_{k=0}^i \frac{(-t)^k}{(1-t)^k} \right )\,dt \\
&= (n+c+1)\left ( \int_0^1 (-1)^i t^{i+d+1}(1-t)^{n+c-d-i} dt + \int_0^1 t^d (1-t)^{n+c+1-d} dt \right ) \\
&= (n+c+1) \left ( (-1)^i \frac{\Gamma(i+d+2)\Gamma(n+c+1-d-i)}{\Gamma(n+c+3)} + \frac{\Gamma(d+1)\Gamma(n+c-d+2)}{\Gamma(n+c+3)} \right ).
\end{align*}
To finish the proof, we apply the relation $B(u,v)=\Gamma(u)\Gamma(v)/\Gamma(u+v)$.
\end{proof}

\begin{remark}
When $i=n$, then Equation \eqref{lem1_id2} becomes Equation (20), Corollary 1, in Batir and Sofo \cite{Batir}.
\end{remark}

\begin{theorem}\label{main_thm2}
For integers $c,d,i$ and $n$ with $0\leq i \leq n$ and $0\leq d \leq c$ and any complex $y$ we have
\begin{equation}\label{main_thm2_id}
\sum_{k=0}^n \frac{(-1)^k}{\binom{n+c}{k+d}} \sum_{j=0}^k \binom{n}{j} y^j = \frac{n+c+1}{n+c+2}\left (\frac{(-1)^n}{\binom{n+c+1}{n+d+1}}
(1+y)^n + \sum_{j=0}^n \frac{\binom{n}{j}}{\binom{n+c+1}{j+d}} (-y)^j \right ).
\end{equation}
\end{theorem}
\begin{proof}
Apply \eqref{Zhang_id} with $f_k=(-1)^k\,\binom{n+c}{k+d}^{-1}$ and $g_k=y^k$. Then, by Lemma \ref{lem1},
\begin{equation*}
S_f(j) = \frac{n+c+1}{n+c+2}\left (\frac{(-1)^{j-1}}{\binom{n+c+1}{j+d}} + \frac{1}{\binom{n+c+1}{d}} \right ).
\end{equation*}
By the direct calculation, this formula also holds for $j=0$.Also,
\begin{equation*}
P(n,0) = \frac{n+c+1}{n+c+2}\left ( \sum_{j=0}^n \binom{n}{j} \frac{(-1)^{j-1}}{\binom{n+c+1}{j+d}} y^j + \frac{(1+y)^n}{\binom{n+c+1}{d}}\right ),
\end{equation*}
and \eqref{main_thm2_id} follows after some obvious simplifications.
\end{proof}

\begin{remark}
Due to the fact that
\begin{equation}\label{GKP_id}
\sum_{j=0}^k \binom{n}{j} (-1)^j = (-1)^k \binom{n-1}{k},
\end{equation}
which is Equation (5.16) in \cite{Graham} or Equation (39) in \cite{Stenlund}, it is clear that when $y=-1$, 
then the double sums in \eqref{main_thm2_id} become single sums involving quotients of binomial coefficients. 
These types of sums have been studied recently by Chu and Guo \cite{ChuGuo}. 
\end{remark}

We proceed with four corollaries that can be deduced from Theorem \ref{main_thm2}. 

\begin{corollary}
For any complex $y$ we have
\begin{align}
&\sum_{k=0}^n (-1)^k \sum_{j=0}^k \frac{\binom{n}{j}}{\binom{n}{k}} y^j \nonumber \\
&\quad = \begin{cases}
  \frac{n+1}{2} + (-1)^n (n+1)\frac{0^n}{n+2}, &\text{if $y=-1$}; \\ 
  \frac{1}{n+2}\left ( (-1)^n (n+1) (1+y)^n + \frac{(-1)^n y^{n+2} + n(1+y)+2y+1}{(1+y)^2}\right ),&\text{if $y\neq -1$}; 
  \end{cases}\label{eq:sum_y^jCor23}
\end{align}
\end{corollary}
\begin{proof}
Set $c=d=0$ in \eqref{main_thm2_id} and calculate
\begin{align*}
\sum_{k=0}^n (-1)^k \sum_{j=0}^k \frac{\binom{n}{j}}{\binom{n}{k}} y^j &= 
\frac{n+1}{n+2} \left ( (-1)^n (1+y)^n + \sum_{j=0}^n \frac{\binom{n}{j}}{\binom{n+1}{j}} (-y)^j \right ) \\
&= \frac{n+1}{n+2} \left ( (-1)^n (1+y)^n + \frac{1}{n+1} \sum_{j=0}^n (n+1-j) (-y)^j \right ).
\end{align*}
When $y\neq -1$, then the last sum is an elementary geometric sum. When $y=-1$, then 
$$\sum_{j=0}^n (n+1-j) = \frac{(n+1)(n+2)}{2}.$$
\end{proof}

As another particular evaluation we offer
$$\sum_{k=0}^n \sum_{j=0}^k \frac{\binom{n}{j}}{\binom{n}{k}} (-1)^{k+j} 2^j = \frac{2}{n+2}\left (2^{n+1} - 1 \right ).$$

\begin{corollary}
For any complex $y$ we have
\begin{align}
&\sum_{k=0}^n (-1)^k \sum_{j=0}^k \frac{\binom{n}{j}}{\binom{n+1}{k}} y^j \nonumber \\
&\quad = \begin{cases}
  \frac{n^2+5n+6}{3(n+3)} + (-1)^n \frac{0^n}{n+3}, &\text{if $y=-1$}; \\ 
  \frac{1}{n+3}\left ( (-1)^n (1+y)^n + \frac{n^2 (1+y)^{2} + n(5y^2+8y+3)+6y^2+6y+2-2(-y)^{n+3}}{(n+1)(1+y)^3}\right ),&\text{if $y\neq -1$}, 
  \end{cases}
\end{align}
\end{corollary}
\begin{proof}
Set $c=1$ and $d=0$ in \eqref{main_thm2_id} and calculate
\begin{align*}
&\sum_{k=0}^n (-1)^k \sum_{j=0}^k \frac{\binom{n}{j}}{\binom{n+1}{k}} y^j = 
\frac{n+2}{n+3} \left ( \frac{(-1)^n}{n+2} (1+y)^n + \sum_{j=0}^n \frac{\binom{n}{j}}{\binom{n+2}{j}} (-y)^j \right ) \\
&\qquad = \frac{n+2}{n+3} \left ( \frac{(-1)^n}{n+2} (1+y)^n + \frac{1}{(n+1)(n+2)} \sum_{j=0}^n (n+1-j)(n+2-j) (-y)^j \right ).
\end{align*}
The sums are again elementary and the stated expression follows after some simplifications.
\end{proof}

Here, we record the double sum evaluation
$$\sum_{k=0}^n \sum_{j=0}^k \frac{\binom{n}{j}}{\binom{n+1}{k}} (-1)^{k+j} 2^j 
= \frac{1}{n+3} \left ( 1 - \frac{n^2 + 7n + 14 - 2^{n+4}}{(n+1)}\right ).$$

\begin{corollary}
If $y$ is a complex number, then
\begin{equation}
\sum_{k=0}^n \frac{(-1)^k}{\binom{2n-1}{k}} \sum_{j=0}^k \binom{n}{j} y^j = \frac{2n}{2n+1}\left (\frac{(-1)^n}{\binom{2n}{n-1}} (1+y)^n 
+ \sum_{j=0}^n \frac{\binom{n}{j}}{\binom{2n}{j}} (-y)^j \right ).
\end{equation}
In particular,
\begin{equation}\label{eq:wqdni23o42}
\sum_{k=0}^n \sum_{j=0}^k \frac{\binom{n}{j}}{\binom{2n-1}{k}} (-1)^{k+j} = \frac{2n}{n+1} \qquad (n\geq 1),
\end{equation}
and
\begin{equation}
\sum_{k=0}^n \sum_{j=0}^k \frac{\binom{n}{j}}{\binom{2n-1}{k}} (-1)^{k+j} 2^j 
= \frac{2n}{2n+1} \left ( \frac{1}{\binom{2n}{n-1}} + \frac{2^{2n}}{\binom{2n}{n}}\right ).
\end{equation}
\end{corollary}
\begin{proof}
Set $c=n-1$ and $d=0$ in \eqref{main_thm2_id}. The first particular case follows from inserting $y=-1$, while the second 
is obtained by inserting $y=-2$ and using
$$\sum_{j=0}^n \frac{\binom{n}{j}}{\binom{2n}{j}} 2^j = \frac{2^{2n}}{\binom{2n}{n}}.$$
\end{proof}

\begin{corollary}
If $y$ is a complex number, then
\begin{equation}
\sum_{k=0}^n \frac{(-1)^k}{\binom{2n-1}{n+k-1}} \sum_{j=0}^k \binom{n}{j} y^j = \frac{2n}{2n+1}\left ((-1)^n(1+y)^n 
+ \sum_{j=0}^n \frac{\binom{n}{j}}{\binom{2n}{j+n-1}} (-y)^j \right ).
\end{equation}
In particular,
\begin{equation}\label{eq:qe2h3798e8dwciebfdqweqw}
\sum_{k=0}^n \sum_{j=0}^k \frac{\binom{n}{j}}{\binom{2n-1}{n+k-1}} (-1)^{k+j} = \frac{2}{n+1} \qquad (n\geq 1),
\end{equation}
\end{corollary}
\begin{proof}
Set $c=d=n-1$ in in \eqref{main_thm2_id}. The particular case is $y=-1$ with the aid of the sum (which is proved below in the lemma)
$$\sum_{j=0}^n\frac{\binom{n}{j}}{\binom{2n}{j+n-1}}=\frac{2n+1}{n(n+1)}.$$
\end{proof}

\begin{lemma}
We have
$$\sum_{j=0}^n \frac{\binom{n}{j}}{\binom{2n}{j+n-1}} = \frac{2n+1}{n(n+1)}.$$
\end{lemma}
\begin{proof}
Since we have
$$
\frac{\binom{n}{j}}{\binom{2n}{n+j-1}} = \frac{(n-j+1)(n+j-1)!}{j!(n+1)\cdots (2n)} = \frac{(n-1)!}{(n+1)\cdots (2n)}(n-j+1)\binom{n+j-1}{j},
$$
the sum follows from the hockey-stick identity for the binomial coefficients:
    \begin{align*}
        \sum_{j=0}^n(n-j+1)\binom{n+j-1}{j} &= \sum_{j=0}^n(n+1)\binom{n+j-1}{j}-\sum_{j=0}^nj\binom{n+j-1}{j}\\
        &= (n+1)\binom{2n}{n}-\sum_{j=1}^nn\binom{n+j-1}{j-1}\\
        &= (n+1)\binom{2n}{n}-n\sum_{j=0}^{n-1}\binom{n+j}{j}\\
        &= (n+1)\binom{2n}{n}-n\binom{2n}{n-1}\\
        &= \binom{2n+1}{n+1}.
    \end{align*}
\end{proof}

\begin{remark}
A closed form evaluation of the sum
$$\sum_{j=0}^n \frac{\binom{n}{j}}{\binom{2n}{j}} (-y)^j$$
does not appear to be possible except for a few particular values of $y$.
\end{remark}

\begin{remark}
The two particular instances of \eqref{main_thm2_id} presented above with $d=0$ and $c\in \{0,1\}$ are special cases of a general 
identity of the similar form. Notice that for $d=0$ and any positive integer $c$, we have
$$\frac{\binom{n}{j}}{\binom{n+c+1}{j}}=\frac{(n-j+1)^{\overline{c+1}}}{(n+1)^{\overline{c+1}}}\quad \text{and}\quad \binom{n+c+1}{n+1}=\frac{(n+2)^{\overline{c}}}{c!},$$
where $x^{\overline{c}}=x\cdot (x+1)\cdots (x+c-1)$ denotes the rising factorial of $x$. With that in mind, we have
\begin{align*}
    \sum_{j=0}^n\frac{\binom{n}{j}}{\binom{n+c+1}{j}}(-y)^j&=\frac{1}{(n+1)^{\overline{c+1}}}\sum_{j=0}^n(n-j+1)^{\overline{c+1}}(-y)^j\\
    &=\frac{1}{(n+1)^{\overline{c+1}}}\sum_{j=1}^{n+1}j^{\overline{c+1}}(-y)^{n-j+1}\\
    &=\frac{(-y)^{n+1}}{(n+1)^{\overline{c+1}}}\sum_{j=1}^{n+1}j^{\overline{c+1}}\left(-\frac{1}{y}\right)^j.
\end{align*}
We can further change the form (for example using Stirling or Eulerian numbers), but let us leave it as is. Thus, we can express 
the case $d=0$ of \eqref{main_thm2_id} as 
\begin{align*}
\sum_{k=0}^n \sum_{j=0}^k  (-1)^k\frac{\binom{n}{j}}{\binom{n+c}{k}}y^j& = \frac{n+c+1}{n+c+2}\left (\frac{(-1)^n}{\binom{n+c+1}{n+1}}
(1+y)^n + \sum_{j=0}^n \frac{\binom{n}{j}}{\binom{n+c+1}{j}} (-y)^j \right )\\
&=\frac{n+c+1}{(n+2)^{\overline{c+1}}}\left((-1)^n(1+y)^n c!+\frac{1}{n+1}\sum_{j=0}^{n+1} j^{\overline{c+1}}(-y)^{n-j+1}\right).
\end{align*}
\end{remark}

\subsection{Some triple sums}

Using some of the identities obtained so far we can deliver triple sums.
\begin{corollary}
    We have
    \begin{equation}\label{eq:wdni23o42v2}
        \sum_{n=1}^m\sum_{k=0}^n \sum_{j=0}^k \frac{\binom{n}{j}}{\binom{2n-1}{k}} (-1)^{k+j} = 2(m+1)-2H_{m+1}
    \end{equation}
    and
    \begin{equation}\label{eq:qe2h3798e8dwciebfdqweqwv2}
\sum_{n=1}^m\sum_{k=0}^n \sum_{j=0}^k \frac{\binom{n}{j}}{\binom{2n-1}{n+k-1}} (-1)^{k+j} = 2H_{m+1}-2,
\end{equation}
\end{corollary}
\begin{proof}
    For \eqref{eq:wdni23o42v2}, sum \eqref{eq:wqdni23o42} for $n=1$ to $m$.
\end{proof}

\begin{remark}
Notice that by adding \eqref{eq:wdni23o42v2} and \eqref{eq:qe2h3798e8dwciebfdqweqwv2} we obtain $2m$. In fact, adding \eqref{eq:wqdni23o42} and \eqref{eq:qe2h3798e8dwciebfdqweqw} we obtain
$$\sum_{k=0}^n\sum_{j=0}^k\binom{n}{j}\left(\frac{1}{\binom{2n-1}{k}}+\frac{1}{\binom{2n-1}{n+k-1}}\right)(-1)^{k+j} = 2.$$
Upon further simplification using \eqref{GKP_id} we obtain an astonishingly simple identity
$$\sum_{k=0}^{n-1} \binom{n-1}{k} \left(\frac{1}{\binom{2n-1}{k}}+\frac{1}{\binom{2n-1}{n+k-1}} \right) = 2.$$
\end{remark}

\section{A non-alternating double sums}

\begin{theorem}\label{thm:shifted-main3} 
Let $n$ and $d$ be non-negative integers, let $y$ be a complex number, and let $a$ be a complex number such that all binomial coefficients appearing below are defined and nonzero. Assume also that $a-d\neq 1$. Then 
\begin{equation}\label{eq:shifted-main3}
\sum_{k=0}^n \frac{1}{\binom{k+a}{k+d}} \sum_{j=0}^k\binom{n}{j}y^j 
= \frac{a-d}{a-d-1} \left( -\frac{(1+y)^n}{\binom{n+a}{n+d+1}} + \sum_{j=0}^n \frac{\binom{n}{j}}{\binom{a+j-1}{j+d}}y^j \right).  
\end{equation} 
In the exceptional case $a=d+1$, the following identity holds:
\begin{equation}\label{eq:harmonic-borderline} 
\sum_{k=0}^n \frac{1}{k+d+1} \sum_{j=0}^k\binom{n}{j}y^j = (H_{n+d+1}-H_d)(1+y)^n - \sum_{j=0}^n \binom{n}{j}(H_{d+j}-H_d)y^j, 
\end{equation} 
with $H_n$ being the $n$th harmonic number.
\end{theorem}
\begin{proof} 
Set \[ f_k=\frac{1}{\binom{k+a}{k+d}}, \qquad g_k=y^k. \] 
Put $q=a-d$. Using the beta integral, we obtain 
\begin{align*} 
S_f(r) &= \sum_{k=0}^{r-1} \frac{1}{\binom{k+a}{k+d}} = q\int_0^1 t^d(1-t)^{q-1} \sum_{k=0}^{r-1}t^k\,dt \\ 
&= q\int_0^1 t^d(1-t)^{q-2}(1-t^r)\,dt = \frac{q}{q-1} \left( \frac{1}{\binom{a-1}{d}} - \frac{1}{\binom{a+r-1}{r+d}} \right). 
\end{align*} 
Hence 
\[ S_f(n+1) = \frac{a-d}{a-d-1} \left( \frac{1}{\binom{a-1}{d}} - \frac{1}{\binom{n+a}{n+d+1}} \right) \] 
and 
\begin{equation*} 
P(n,0) = \sum_{j=0}^n\binom{n}{j}y^j S_f(j) = \frac{a-d}{a-d-1} \left( \frac{(1+y)^n}{\binom{a-1}{d}} - \sum_{j=0}^n \frac{\binom{n}{j}}{\binom{a+j-1}{j+d}}y^j \right). 
\end{equation*} 
Substitution into Zhang's identity yields \eqref{eq:shifted-main3}, because the terms containing $\binom{a-1}{d}^{-1}(1+y)^n$ cancel. 

It remains to consider the exceptional case $a=d+1$. In this case, we have
\[f_k=\frac{1}{k+d+1},\qquad S_f(r)=\sum_{k=0}^{r-1}\frac{1}{k+d+1}=H_{d+r}-H_d.\]
Furthermore,
\[H_g(n)=(1+y)^n,\qquad P(n,0)=\sum_{j=0}^n\binom{n}{j}(H_{d+j}-H_d)y^j.\]
An application of Zhang's identity now gives \
\begin{align*} 
\sum_{k=0}^n \frac{1}{k+d+1} \sum_{j=0}^k\binom{n}{j}y^j &= (H_{n+d+1}-H_d)(1+y)^n- \sum_{j=0}^n \binom{n}{j}(H_{d+j}-H_d)y^j,
\end{align*} 
which is \eqref{eq:harmonic-borderline}.
\end{proof}

In view of an identity in Gould's book~\cite[Identity (4.13)]{Gould}, a full evaluation of the double sum in~\eqref{eq:shifted-main3} is presented in Theorem~\ref{thm.veuifox} for the case when $d=0$ and $a=m+1$ is a positive integer.
\begin{theorem}\label{thm.veuifox}
If $n$ is a non-negative integer, $m$ is a positive integer and $y$ is a complex number, then
\begin{align}
&\sum_{k = 0}^n {\frac{1}{{\binom{{k + 1 + m}}{k}}}\sum_{j = 0}^k {\binom{{n}}{j}y^{j + m} } } \nonumber\\
&\qquad = \frac{{m + 1}}{m}\left( {y^m  + \frac{{\left( {1 + y} \right)^n }}{{\binom{{n + m}}{n}}}\left( {\left( {1 + y} \right)^m  - \frac{{n + 1}}{{n + 1 + m}}y^m } \right) - \frac{1}{{\binom{{n + m}}{m}}}\sum_{k = 0}^m {\binom{{n + m}}{k}y^k } } \right).
\end{align}
\end{theorem}

\begin{corollary}
We have
\begin{equation}\label{cor_thm3_id1}
\sum_{k=0}^n \frac{\binom{n}{k}}{\binom{k+1+c}{k}} (-1)^k = \frac{c+1}{c+1+n}.
\end{equation}
\end{corollary}
\begin{proof}
Work with 
\begin{equation}\label{main_thm3_id}
\sum_{k=0}^n \frac{1}{\binom{k+1+c}{k}} \sum_{j=0}^k \binom{n}{j} y^j = \frac{c+1}{c}\left ( - \frac{(1+y)^n}{\binom{n+c+1}{n+1}} 
+ \sum_{j=0}^n \frac{\binom{n}{j}}{\binom{c+j}{j}} y^j \right ),
\end{equation}
which is a special case of \eqref{eq:shifted-main3} for $a=c+1$ and $d=0$. Set $y=-1$ in \eqref{main_thm3_id}, use \eqref{GKP_id}, 
and shift the summation index from $n-1$ to $n$.
\end{proof}

\begin{corollary}
We have
\begin{equation}\label{cor_thm3_id2}
\sum_{k=0}^n \sum_{j=0}^k \frac{\binom{n}{j}}{(k+1)(k+2)} = 2^n \frac{n+3}{(n+1)(n+2)} - \frac{1}{n+1}.
\end{equation}
\end{corollary}
\begin{proof}
Set $y=c=1$ in \eqref{main_thm3_id}.
\end{proof}

\begin{corollary}
We have
\begin{equation}\label{cor_thm3_id3}
\sum_{k=0}^n \sum_{j=0}^k \frac{\binom{n}{j}}{\binom{n+1+k}{k}} 
= \frac{n+1}{n}\left ( \frac{2^{2n-1}}{\binom{2n}{n}} - \frac{n+1}{2n+1}\frac{2^n}{\binom{2n}{n}} + \frac{1}{2} \right ).
\end{equation}
\end{corollary}
\begin{proof}
Set $y=1$ and $c=n$ in \eqref{main_thm3_id} and use the evaluation
$$ \sum_{j=0}^n \frac{\binom{n}{j}}{\binom{n+j}{j}} = \frac{2^{2n-1}}{\binom{2n}{n}} + \frac{1}{2}.$$
\end{proof}

\begin{proposition}
If $n$ is a non-negative integer and $y$ is a complex number, then
\begin{equation}\label{to5m7l9}
\sum_{k = 0}^n 2^{2k} \binom{{2k}}{k}^{ - 1} \sum_{j = 0}^k \binom{{n}}{j} y^j = \frac{1}{3} + \frac{{2^{2n + 1} }}{{3C_n }}\left( {1 + y} \right)^n - \frac{1}{6}\sum_{k = 1}^n {\frac{{2^{2k} }}{{C_{k - 1} }}\binom{{n}}{k}y^k }.
\end{equation}
\end{proposition}
\begin{proof}
Set $c=-3/2$ in~\eqref{main_thm3_id} and evaluate the resulting generalized binomial coefficients using the Gamma function.
\end{proof}

Let $(s_n)_{n\geq 0}$ and $(\sigma_n)_{n\geq 0}$ be two sequences of complex numbers. The pair $\{(s_n),(\sigma_n) \}$ is a binomial transform pair if these sequences are connected by the relations
\begin{equation*}
s_n = \sum_{k = 0}^n \binom{n}{k} (-1)^k \sigma_k \iff \sigma_n = \sum_{k = 0}^n \binom{n}{k} (-1)^k s_k.
\end{equation*}
\begin{lemma}[{\cite[Lemma 2.1]{Adegoke0}}]\label{qqdca04}
Let $\left\{(s_n),(\sigma_n)\right\}$, $n=0,1,2,\ldots$, be a binomial-transform pair. Let $\mathcal L_x$ be a linear operator defined by $\mathcal L_x(x^j)=s_j$ for every complex number $x$ and every non-negative integer~$j$. Then $\mathcal L_x((1-x)^j)=\sigma_j$.

In particular, if  $\{(s_n), (\sigma_n)\}$, $n=0,1,2,\ldots$, is a binomial-transform pair, then so is $\{(a_n), (\alpha_n)\}$,  where
\begin{equation}\label{w95ejl1}
a_n  = s_{n + m} \text{ and }\alpha _n  = \sum_{p = 0}^m {( - 1)^p \binom{{m}}{p}\sigma _{n + p} }. 
\end{equation}
\end{lemma}

\begin{proposition}
Let $\{(s_n),(\sigma_n)\}$, $n=0,1,2,\dots$, be a binomial-transform pair. If $m$ is a non-negative integer, then
\begin{align}\label{v5beyse}
&\sum_{k = 0}^n {2^{2k} \binom{{2k}}{k}^{ - 1} \sum_{j = 0}^k {(-1)^j\binom{{n}}{j}s_{m + j} } } \nonumber\\
&\qquad = \frac{1}{3}s_m  + \frac{{2^{2n + 1} }}{{3C_n }}\sum_{k = 0}^m {( - 1)^k \binom{{m}}{k}\sigma _{k + n} }  + \frac{1}{6}\sum_{k = 1}^n {(-1)^{k-1}\frac{{2^{2k} }}{{C_{k - 1} }}\binom{{n}}{k}s_{m + k} } .
\end{align}
In particular,
\begin{equation}
\sum_{k = 0}^n {2^{2k} \binom{{2k}}{k}^{ - 1} \sum_{j = 0}^k {(-1)^j\binom{{n}}{j}s_j } }  = \frac{1}{3}s_0  + \frac{{2^{2n + 1} }}{{3C_n }}\sigma _n  + \frac{1}{6}\sum_{k = 1}^n {( - 1)^{k - 1} \frac{{2^{2k} }}{{C_{k - 1} }}\binom{{n}}{k}s_k } .\label{sv9ff0x}
\end{equation}
\end{proposition}
\begin{proof}
Write $-x$ for $y$ in~\eqref{to5m7l9}, mutliply both sides by $x^m$ and act on both sides with $\mathcal L_x$. This gives
\begin{align*} 
\sum_{k=0}^n \frac{2^{2k}}{\binom{2k}{k}} \sum_{j=0}^k (-1)^j\binom{n}{j} \mathcal L_x(x^{m+j})&= \frac13\mathcal L_x(x^m) + \frac{2^{2n+1}}{3C_n} \mathcal L_x\bigl(x^m(1-x)^n\bigr)\\
&\qquad+ \frac16\sum_{k=1}^n (-1)^{k-1} \frac{2^{2k}}{C_{k-1}} \binom{n}{k}\mathcal L_x(x^{m+k}). 
\end{align*}
Furthermore, $\mathcal L_x(x^{m+j})=s_{m+j}$ and, by Lemma~\ref{qqdca04}, \[ \mathcal L_x\bigl(x^m(1-x)^n\bigr) = \sum_{p=0}^m (-1)^p\binom{m}{p}\sigma_{n+p}. \] 
Substitution gives \eqref{v5beyse}. The case $m=0$ gives \eqref{sv9ff0x}.
\end{proof}
A binomial transform version of Identity \eqref{main_thm3_id} is given in the next theorem.

\begin{theorem}
Let $\{(s_k),(\sigma_k)\}$, $k=0,1,2,\dots$, be a binomial-transform pair. If $n$ is a non-negative integer and $c$ is a complex number that is not a non-positive integer, then
\begin{equation}\label{ml097wi}
\sum_{k=0}^n \frac{1}{\binom{k+1+c}{k}} \sum_{j=0}^k (-1)^j\binom{n}{j} s_j = \frac{c+1}{c}\left ( - \frac{\sigma_n}{\binom{n+c+1}{n+1}} 
+ \sum_{k=0}^n (-1)^k\frac{\binom{n}{k}}{\binom{c+k}{k}} s_k \right ).
\end{equation}
\end{theorem}
\begin{proof}
    Replace $y$ by $-x$ in \eqref{main_thm3_id}. Applying $\mathcal L_x$ to both sides and use Lemma~\ref{qqdca04}. We omit the details.
\end{proof}

Many new identities involving special numbers can be derived from~\eqref{v5beyse} and~\eqref{ml097wi}. We illustrate with a couple of identities based on~\eqref{sv9ff0x}.
\begin{proposition}\label{wnwq8ha}
If $n$ is a non-negative integer, then
\begin{align}
\sum_{k = 0}^{2n} {2^{2k} \binom{{2k}}{k}^{ - 1} \sum_{j = 1}^k {( - 1)^j \binom{{2n}}{j}2^j } }  &= \frac{1}{6}\sum_{k = 1}^{2n} {( - 1)^{k - 1} \frac{{2^{3k} }}{{C_{k - 1} }}\binom{{2n}}{k}}\label{vqui2ye} ,\\
\sum_{k = 0}^{2n} {2^{2k} \binom{{2k}}{k}^{ - 1} \sum_{j = 1}^k {( - 1)^j \binom{{2n}}{j}\frac{{2^j H_j }}{{j + 1}}} }  &= \frac{1}{6}\sum_{k = 1}^{2n} {( - 1)^{k - 1} \frac{{2^{3k} }}{{C_{k - 1} }}\binom{{2n}}{k}\frac{{H_k }}{{k + 1}}}\label{axngxm8} ,
\end{align}
and
\begin{equation}\label{r0mtltv}
\sum_{k = 0}^{2n} {2^{2k} \binom{{2k}}{k}^{ - 1} \sum_{j = 1}^k {( - 1)^j \binom{{2n}}{j}j\,2^{ - j} C_j } }  = \frac{1}{3}\sum_{k = 1}^{2n} {( - 1)^{k - 1} 2^k \binom{{2n}}{k}\frac{{k\left( {2k - 1} \right)}}{{k + 1}}} .
\end{equation}
\end{proposition}
\begin{proof}
From the elementary identity
\begin{equation*}
\sum_{k = 1}^n {( - 1)^k \binom{{n}}{k}2^k }  = ( - 1)^n  - 1,
\end{equation*}
we identify the binomial transform pair $\{(s_n),(\sigma_n)\}$, where
\begin{equation}\label{riiuf2p}
s_0=0,\quad s_n=2^n,\,n>0,\quad \sigma_{2n}=0,\text{ and }\sigma_{2n+1}=-2.
\end{equation}
Use~\eqref{riiuf2p} in~\eqref{sv9ff0x} after replacing $n$ with $2n$ to obtain~\eqref{vqui2ye}. To prove~\eqref{axngxm8}, engage the binomial transform pair (see~\cite{rathie22}):
\begin{equation*}
s_n = \frac{2^nH_n}{n+1} \quad\text{and}\quad \sigma_n = -\frac{1-(-1)^n}{\,\left(n+1\right)}\,O_{(n+1)/2}.
\end{equation*}
Finally, to derive~\eqref{r0mtltv}, replace $n$ with $2n$ in~\eqref{sv9ff0x} and use the binomial transform pair (see~\cite{sulei23}):
\begin{equation*}
s_n = n2^{- n} C_{n} ,\quad\text{and}\quad \sigma_n 
= \begin{cases}
 -n2^{- n} C_{(n-1)/2},&\text{if $n$ is odd;}  \\ 
 0,&\text{if $n$ is even.}  
 \end{cases}
\end{equation*}
\end{proof}

\begin{proposition}
If $n$ is a non-negative integer, $c\in\mathbb C\setminus\{0,-1,-2,\ldots\}$, and $y$ is a complex number, then
\begin{align}\label{vxhdwmb}
\sum_{k = 0}^n {\frac{{H_{k + c + 1} }}{{\binom{{k + c + 1}}{k}}}\sum_{j = 0}^k {\binom{{n}}{j}y^j } }  = \frac{{c + 1}}{c}\left( {H_c  + \frac{1}{c}} \right) - \frac{{c + 1}}{c}\frac{{H_{n + c + 1}  + 1/c}}{{\binom{{n + c + 1}}{c}}}\left( {1 + y} \right)^n \nonumber\\
 + \frac{{c + 1}}{c}\sum_{k = 1}^n {\frac{{H_{k + c}  + 1/c}}{{\binom{{c + k}}{k}}}\binom{{n}}{k}y^k } .
\end{align}
\end{proposition}
\begin{proof}
Differentiate~\eqref{main_thm3_id} with respect to $c$ to obtain
\begin{align}\label{a2yt3mz}
\sum_{k = 0}^n {\frac{{\left( {H_{k + c + 1}  - H_{c + 1} } \right)}}{{\binom{{k + c + 1}}{k}}}\sum_{j = 0}^k {\binom{{n}}{j}y^j } }  &= \frac{1}{{c^2 }} - \frac{{\left( {1 + y} \right)^n }}{{\binom{{n + c + 1}}{c}}}\left( {\frac{1}{{c^2 }} + \frac{{c + 1}}{c}\left( {H_{n + c + 1}  - H_c } \right)} \right)\nonumber\\
&\qquad + \sum_{k = 1}^n {\left( {\frac{1}{{c^2 }} + \frac{{c + 1}}{c}\left( {H_{k + c}  - H_c } \right)} \right)\frac{{\binom{{n}}{k}}}{{\binom{{c + k}}{k}}}y^k }, 
\end{align}
from which~\eqref{vxhdwmb} follows after using~\eqref{main_thm3_id} again.
\end{proof}

\begin{proposition}
If $n$ is a non-negative integer and $y$ is a complex number, then
\begin{align}
\sum_{k = 0}^n {2^{2k} \binom{{2k}}{k}^{ - 1} O_k \sum_{j = 0}^k {\binom{{n}}{j}y^j} }  &= \frac{2}{9} + \frac{{2^{2n + 1} }}{{3C_n }}\left( {O_n  - \frac{1}{3}} \right)\left( {1 + y} \right)^n\nonumber\\
&\qquad  - \frac{1}{3}\sum_{k = 1}^n {\frac{{2^{2k - 1} }}{{C_{k - 1} }}\left( {O_{k - 1}  - \frac{1}{3}} \right)\binom{{n}}{k}y^k } .\label{uirpuki}
\end{align}
\end{proposition}
\begin{proof}
Set $c=-3/2$ in~\eqref{a2yt3mz}.
\end{proof}

Corresponding to identity~\eqref{uirpuki}, we have the following binomial transform identity.
\begin{proposition}\label{ch1w2d1}
Let $\{(s_n),(\sigma_n)\}$ be a binomial-transform pair. Then
\begin{align}
\sum_{k = 0}^n {2^{2k} \binom{{2k}}{k}^{ - 1} O_k \sum_{j = 0}^k {(-1)^j\binom{{n}}{j}s_j} }  &= \frac{2}{9}s_0 + \frac{{2^{2n + 1} }}{{3C_n }}\left( {O_n  - \frac{1}{3}} \right)\sigma_n\nonumber\\
&\qquad  + \frac{1}{3}\sum_{k = 1}^n {(-1)^{k-1}\frac{{2^{2k - 1} }}{{C_{k - 1} }}\left( {O_{k - 1}  - \frac{1}{3}} \right)\binom{{n}}{k}s_k } .
\end{align}
\end{proposition}
We close this section with some examples derived from Proposition~\ref{ch1w2d1}. The proofs are analogous to those of Proposition~\ref{wnwq8ha}.
\begin{proposition}
If $n$ is a non-negative integer, then
\begin{align}
\sum_{k = 0}^{2n} {2^{2k} \binom{{2k}}{k}^{ - 1} O_k \sum_{j = 1}^k {( - 1)^j \binom{{2n}}{j}2^j } }  &= \frac{1}{6}\sum_{k = 1}^{2n} {( - 1)^{k - 1} \binom{{2n}}{k}\frac{{2^{3k} }}{{C_{k - 1} }}\left( {O_{k - 1}  - \frac{1}{3}} \right)} ,\\
\sum_{k = 0}^{2n} {2^{2k} \binom{{2k}}{k}^{ - 1} O_k \sum_{j = 1}^k {( - 1)^j \binom{{2n}}{j}\frac{{2^j H_j }}{{j + 1}}} }  &= \frac{1}{6}\sum_{k = 1}^{2n} {( - 1)^{k - 1} \frac{{2^{3k} }}{{C_{k - 1} }}\binom{{2n}}{k}\left( {O_{k - 1}  - \frac{1}{3}} \right)\frac{{H_k }}{{k + 1}}} ,
\end{align}
and
\begin{equation}
\sum_{k = 0}^{2n} {2^{2k} \binom{{2k}}{k}^{ - 1} O_k \sum_{j = 1}^k {( - 1)^j \binom{{2n}}{j}j\,2^{ - j} C_j } }  = \frac{1}{3}\sum_{k = 1}^{2n} {( - 1)^{k - 1} 2^k \binom{{2n}}{k}\frac{{k\left( {2k - 1} \right)}}{{k + 1}}\left( {O_{k - 1}  - \frac{1}{3}} \right)} .
\end{equation}
\end{proposition}

In the final part of this section, recall that Bat\i r and Sofo \cite{Batir} obtained the following identity:
\begin{equation}\label{eq:BaSo2021}
    \sum_{k=0}^n\frac{1}{\binom{n+s}{k}}=\frac{n+s+1}{2^{n+1}}\left(\sum_{k=0}^n\frac{2^k}{(k+1)\binom{s+k+1}{k+1}}+\sum_{k=0}^n\frac{2^k}{k+s+1}\right).
\end{equation}
\begin{proposition}
    We have
    \begin{align*}
        \sum_{k=0}^n\sum_{j=0}^k\frac{\binom{n}{j}}{\binom{n+s}{k}}y^j&=\frac{n+s+1}{2^{n+1}}\left(\sum_{k=0}^n\frac{2^k}{(k+1)\binom{s+k+1}{k+1}}+\sum_{k=0}^n\frac{2^k}{k+s+1}\right)(1+y)^n\\
        &\qquad -\sum_{j=0}^n\binom{n}{j}y^j\sum_{k=0}^{j-1}\frac{1}{\binom{n+s}{k}}
    \end{align*}
\end{proposition}
\begin{proof}
    Set $f_k=\binom{n+s}{k}^{-1}$ and $g_k=y^k$ in \eqref{Zhang_id}, which gives
    $$\sum_{k=0}^n\sum_{j=0}^k\frac{\binom{n}{j}}{\binom{n+s}{k}}y^j=(1+y)^n\sum_{k=0}^n\frac{1}{\binom{n+s}{k}}-\sum_{j=0}^n\binom{n}{j}y^j\sum_{k=0}^{j-1}\frac{1}{\binom{n+s}{k}}.$$
    Make use of \eqref{eq:BaSo2021}. 
    \end{proof}

\begin{remark}
In general, this partial sum does not appear to admit a simplification comparable to  \eqref{eq:BaSo2021}:
$$\sum_{k=0}^{j-1}\frac{1}{\binom{n+s}{k}}.$$
The sum 
$$\sum_{j=0}^n\binom{n}{j}y^j\sum_{k=0}^{j-1}\frac{1}{\binom{n+s}{k}}$$
can be changed by swapping the order of summation to
$$\sum_{j=0}^n \binom{n}{j} y^j \sum_{k=0}^{j-1} \frac{1}{\binom{n+s}{k}} 
= \sum_{k=0}^{n-1} \frac{1}{\binom{n+s}{k}} \sum_{j=k+1}^n \binom{n}{j} y^j.$$
\end{remark}

\section{Quadratic incomplete sums and Legendre polynomials}

In Zhang's paper \cite{Zhang1}, the following identity for squared incomplete binomial sums is established: 
\begin{equation}\label{eq:Zhang-quadratic} 
\sum_{k=0}^n f_k \left( \sum_{j=0}^k\binom{n}{j}g_j \right)^2 = S_f(n+1)H_g(n)^2 - H_g(n)P(n,0) - \sum_{k=0}^n\binom{n}{k}\Delta(n,k), 
\end{equation} 
where 
\begin{equation*}
S_f(r) = \sum_{k=0}^{r-1}f_k, \qquad H_g(n) = \sum_{j=0}^n\binom{n}{j}g_j, \qquad P(n,k) = \sum_{j=k}^n\binom{n}{j}g_jS_f(j),   
\end{equation*}
and
\begin{equation*}
\Delta(n,k) = g_kP(n,k)-g_{n-k}Q(n,k), \qquad Q(n,k) = \sum_{j=k}^n \binom{n}{j}g_{n-j}S_f(n-j).
\end{equation*}
This identity received some interest recently and was further generalized by the third author in \cite{Gryszka2026}.

We can apply \eqref{eq:Zhang-quadratic} to the sequences 
\[ f_k=\frac{(-1)^k}{\binom{n+c}{k+d}}, \qquad g_k=(-1)^k, \] 
where $n\geq 1$ and $0\leq d\leq c$. Making use again of \eqref{GKP_id} we obtain the equivalent form 
\begin{equation}\label{eq:quadratic-reduced-sum} 
\sum_{k=0}^n \frac{(-1)^k}{\binom{n+c}{k+d}} \left( \sum_{j=0}^k (-1)^j \binom{n}{j} \right)^2 
= \sum_{k=0}^{n-1} (-1)^k \frac{\binom{n-1}{k}^2}{\binom{n+c}{k+d}}. 
\end{equation} 

Our main result in this section is an integral representation of the sum on the right in terms of the Legendre polynomials. We recall only the properties needed below; for a comprehensive treatment, see \cite[Chapters 14 and 18]{NISTDLMF}, \cite{OlverEtAl2010}, or the classical monograph \cite{Szego1975}.

The Legendre polynomials $P_n(x)$ are defined as an orthogonal system with respect to the weight function $w(x)=1$ 
over the interval $[-1,1]$. Rodrigues' formula for $P_n(x)$ is 
\[ P_n(x) = \frac{1}{2^n n!}\frac{d^n}{dx^n}(x^2-1)^n \]
and one of many equivalent expressions is
\begin{equation}\label{eq:Legendre-explicit}
P_n(x) = \frac{1}{2^n}\sum_{k=0}^n \binom{n}{k}^2 (x-1)^{n-k} (x+1)^k.
\end{equation}

\begin{theorem}\label{thm:quadratic-Legendre} 
Let $n$ be a positive integer and let $c$ and $d$ be such that $0\leq d\leq c$. Then
\begin{equation}\label{eq:quadratic-Legendre}
\sum_{k=0}^{n-1} (-1)^k \frac{\binom{n-1}{k}^2}{\binom{n+c}{k+d}} = (n+c+1) \int_0^1 t^d(1-t)^{c-d+1} P_{n-1}(1-2t)\,dt,  
\end{equation} 
where $P_n(x)$ denotes the Legendre polynomial of degree $n$. 
\end{theorem} 
\begin{proof} 
Using the Beta integral representation, we have 
\begin{align*} 
\sum_{k=0}^{n-1} (-1)^k \frac{\binom{n-1}{k}^2}{\binom{n+c}{k+d}} 
&= (n+c+1) \sum_{k=0}^{n-1} (-1)^k\binom{n-1}{k}^2 \int_0^1 t^{k+d}(1-t)^{n+c-k-d}\,dt\\ 
&= (n+c+1) \int_0^1 t^d(1-t)^{c-d+1} \sum_{k=0}^{n-1} (-1)^k\binom{n-1}{k}^2 t^k(1-t)^{n-1-k}\,dt. 
\end{align*} 
We now set $x=1-2t$ in \eqref{eq:Legendre-explicit} and obtain the shifted identity
\begin{equation}\label{eq:shifted-Legendre-binomial} 
P_m (1-2t) = \sum_{k=0}^m (-1)^k \binom{m}{k}^2 t^k (1-t)^{m-k}. 
\end{equation} 
Taking $m=n-1$ in \eqref{eq:shifted-Legendre-binomial} gives \eqref{eq:quadratic-Legendre}. 
\end{proof} 
\begin{remark}
    The integral representation \eqref{eq:quadratic-Legendre} extends to complex parameters $c$ and $d$ under suitable convergence conditions, for example 
    \[ \Re d>-1, \qquad \Re(c-d)>-2, \] 
    provided that the generalized binomial coefficients involved are defined.
\end{remark}

The integral representation allows us to establish a large family of vanishing sums. 
\begin{corollary}\label{cor:quadratic-vanishing} 
Let $c,d$ be integers satisfying $0\leq d\leq c$. If $n\geq c+3,$ then 
\begin{equation}\label{eq:quadratic-vanishing} 
\sum_{k=0}^n \frac{(-1)^k}{\binom{n+c}{k+d}} \left( \sum_{j=0}^k (-1)^j \binom{n}{j} \right)^2 
= \sum_{k=0}^{n-1} (-1)^k \frac{\binom{n-1}{k}^2}{\binom{n+c}{k+d}} = 0.
\end{equation} 
In particular, for $n\geq 3$
\begin{equation}
\sum_{k=0}^n \frac{(-1)^k} {\binom{n}{k}} \left( \sum_{j=0}^k(-1)^j\binom{n}{j} \right)^2 
= \sum_{k=0}^{n-1} (-1)^k \frac{\binom{n-1}{k}^2}{\binom{n}{k}} = 0. 
\end{equation}
\end{corollary} 
\begin{proof} 
The polynomial 
\[ t^d(1-t)^{c-d+1} \] 
has degree $c+1$. If $n\geq c+3$, then 
\[ c+1\leq n-2. \] 
The shifted Legendre polynomial $P_{n-1}(1-2t)$ is orthogonal on $[0,1]$ to every polynomial of degree at most $n-2$. Consequently, 
\begin{equation*}
\int_0^1 t^d(1-t)^{c-d+1} P_{n-1}(1-2t)\,dt =0. 
\end{equation*}
Identity \eqref{eq:quadratic-vanishing} now follows from Theorem~\ref{thm:quadratic-Legendre} and identity \eqref{eq:quadratic-reduced-sum}. 
\end{proof} 

We can reformulate identity~\eqref{eq:quadratic-vanishing} as a vanishing sum of a slightly different form.

\begin{corollary}\label{cor:quadratic-companion-vanishing} 
Let $c,d$ be integers satisfying $0\leq d\leq c$. If $n\geq c+3$, then 
\begin{equation}\label{eq:quadratic-companion-vanishing} 
\sum_{k=0}^n (-1)^k (n-2k) \frac{\binom{n}{k}^2}{\binom{n+c+1}{k+d}} =0.
\end{equation} 
\end{corollary}
\begin{proof} 
We first establish the polynomial identity \begin{equation}\label{eq:Legendre-companion-polynomial} 
\sum_{k=0}^n (-1)^k(n-2k)\binom{n}{k}^2 t^k(1-t)^{n-k} = nP_{n-1}(1-2t). 
\end{equation} 
Indeed, upon noticing that 
\[ (n-2k)\binom{n}{k} = n\left( \binom{n-1}{k}-\binom{n-1}{k-1} \right), \] 
we obtain 
\begin{align*} 
&\sum_{k=0}^n (-1)^k(n-2k)\binom{n}{k}^2 t^k(1-t)^{n-k}\\ 
&\quad= n\sum_{k=0}^n (-1)^k \left( \binom{n-1}{k}-\binom{n-1}{k-1} \right) \binom{n}{k}t^k(1-t)^{n-k}\\
&\quad=n\sum_{k=0}^n(-1)^k\left(\binom{n-1}{k}^2-\binom{n-1}{k-1}^2\right)t^k(1-t)^{n-k}\\
&\quad=n(1-t)\sum_{k=0}^{n-1}(-1)^k\binom{n-1}{k}^2t^k(1-t)^{n-1-k}\\
&\qquad +nt\sum_{k=0}^{n-1}(-1)^k\binom{n-1}{k}^2t^k(1-t)^{n-1-k}\\
&\quad=n\big((1-t)+t\big)\sum_{k=0}^{n-1}(-1)^k\binom{n-1}{k}^2t^k(1-t)^{n-1-k}\\
&\quad=nP_{n-1}(1-2t).
\end{align*} 
The last equality follows from \eqref{eq:shifted-Legendre-binomial}. This proves \eqref{eq:Legendre-companion-polynomial}.

We apply the beta integral and \eqref{eq:Legendre-companion-polynomial} to obtain
\begin{align}
    \nonumber&\sum_{k=0}^n(-1)^k(n-2k)\frac{\binom{n}{k}^2}{\binom{n+c+1}{k+d}}\\
    \nonumber&\quad=(n+c+2)\int_0^1t^d(1-t)^{c-d+1}\cdot \sum_{k=0}^n(-1)^k(n-2k)\binom{n}{k}^2t^k(1-t)^{n-k}dt\\
    &\quad=n(n+c+2)\int_0^1t^d(1-t)^{c-d+1}P_{n-1}(1-2t)dt. \label{eq:quadratic-companion-integral}
\end{align}
On the other hand, Theorem~\ref{thm:quadratic-Legendre} implies that 
\begin{align*} 
\sum_{k=0}^{n-1} (-1)^k \frac{\binom{n-1}{k}^2}{\binom{n+c}{k+d}}= (n+c+1) \int_0^1 t^d(1-t)^{c-d+1} P_{n-1}(1-2t)\,dt. \end{align*} 
Comparing this identity with \eqref{eq:quadratic-companion-integral}, we obtain 
\begin{align} 
\sum_{k=0}^n (-1)^k(n-2k) \frac{\binom{n}{k}^2}{\binom{n+c+1}{k+d}}= \frac{n(n+c+2)}{n+c+1} \sum_{k=0}^{n-1} (-1)^k \frac{\binom{n-1}{k}^2}{\binom{n+c}{k+d}}. \label{eq:quadratic-companion-comparison} 
\end{align} 
Therefore if $n\geq c+3$, then the sum on the right-hand side of \eqref{eq:quadratic-companion-comparison} vanishes by \eqref{eq:quadratic-vanishing}. Hence the sum on the left-hand side also vanishes, proving \eqref{eq:quadratic-companion-vanishing}.
\end{proof}

Starting from~\eqref{eq:quadratic-Legendre} and using a different representation of shifted Legendre polynomials, namely,
\begin{equation}\label{wc6p34s}
P_n \left( {1 - 2t} \right) = \sum_{k = 0}^n {( - 1)^k \binom{{n}}{k}\binom{{n + k}}{k}t^k },
\end{equation}
we can derive a general result that is valid for a larger family of complex numbers $c$ and $d$.
\begin{corollary}
If $n$ is a non-negative integer and $c$ and $d$ are complex numbers such that $c-1-d$ is not a negative integer and $c\neq d$, then
\begin{equation}\label{j3ybkdj}
\sum_{k = 0}^n {( - 1)^k \frac{{\binom{{n}}{k}^2 }}{{\binom{{n + c - 1}}{{k + d}}}}}  = \frac{{c + n}}{{c - d}}\,\sum_{k = 0}^n {( - 1)^k \frac{{\binom{{n}}{k}\binom{{n + k}}{k}}}{{\binom{{k + c}}{{k + d}}}}} .
\end{equation}
\end{corollary}
\begin{proof}
Substitute~\eqref{wc6p34s} into~\eqref{eq:quadratic-Legendre} and perform termwise integration using Euler's Beta function. When done, replace $n$ with $n+1$ and $c$ with $c-2$ everywhere.
\end{proof}

One can derive some interesting special cases of~\eqref{j3ybkdj} by carefully choosing $c$ and $d$. One can also derive identities involving harmonic numbers by differentiating with respect to $c$ and $d$.

\begin{corollary}
If $n$, $c$, and $d$ are non-negative integers such that $n\ge c$ and $d\le c-2$, then
\begin{equation}
\sum_{k = 0}^n {( - 1)^k \frac{{\binom{{n}}{k}\binom{{n + k}}{k}}}{{\binom{{k + c}}{{k + d}}}}}  = 0.
\end{equation}
\end{corollary}
\begin{proof}
A consequence of~\eqref{eq:quadratic-vanishing}  and~\eqref{j3ybkdj}.
\end{proof}

The integral involving Legendre polynomials can be further computed with the aid of the following Lemma.
\begin{lemma}\label{lem:Legendre-moment} 
Let $m$ and $q$ be a non-negative integer. Then 
\begin{equation}\label{eq:Legendre-moment} 
\int_0^1t^qP_m(1-2t)\,dt = 
\begin{cases} 0, & q\in\{0,1,\ldots,m-1\},\\[2mm] 
\displaystyle (-1)^m \frac{(q!)^2}{(q-m)!(q+m+1)!}, & otherwise. 
\end{cases} 
\end{equation} 
More generally, if $q\in\mathbb C$ satisfies $\Re q>-1$ and $q-m+1$ is not a pole of the Gamma function, then 
\begin{equation}
\label{eq:Legendre-moment-complex} 
\int_0^1t^qP_m(1-2t)\,dt = (-1)^m \frac{\Gamma(q+1)^2} {\Gamma(q-m+1)\Gamma(q+m+2)}. 
\end{equation}
\end{lemma} 
\begin{proof} 
For $q\in\mathbb C$, with $\Re q>-1$, integrating by parts $m$ times and application of 
\[P_m(1-2t)=\frac{1}{m!}\frac{d^m}{dt^m}(t^m(1-t)^m)\]
gives 
\begin{align*} 
\int_0^1t^qP_m(1-2t)dt &= \frac{(-1)^m}{m!} \frac{\Gamma(q+1)}{\Gamma(q-m+1)} \int_0^1t^q(1-t)^mdt\\ 
&= \frac{(-1)^m}{m!} \frac{\Gamma(q+1)}{\Gamma(q-m+1)} B(q+1,m+1)\\ 
&= (-1)^m \frac{\Gamma(q+1)^2} {\Gamma(q-m+1)\Gamma(q+m+2)}. 
\end{align*} 
For $q=0,1,\ldots,m-1$, the expression vanishes either by orthogonality argument or by interpreting $1/\Gamma(q-m+1)$ as zero.
\end{proof} 

As a consequence, we have the following evaluation.
\begin{theorem}\label{thm:quadratic-finite-evaluation} 
Let $n,c,d$ be integers satisfying $n\geq1$ and  $0\leq d\leq c$.Then 
\begin{align} 
&\sum_{k=0}^{n-1} (-1)^k \frac{\binom{n-1}{k}^2}{\binom{n+c}{k+d}} \nonumber\\ 
&\quad= (n+c+1)(-1)^{n-1} \sum_{j=\max\{0,n-1-d\}}^{c-d+1} (-1)^j\binom{c-d+1}{j} \frac{((d+j)!)^2} {(d+j-n+1)!(d+j+n)!}. 
\label{eq:quadratic-finite-evaluation} 
\end{align} 
An empty sum on the right-hand side is understood to be zero. 
\end{theorem} 
\begin{proof} 
Expand 
\[ t^d(1-t)^{c-d+1} = \sum_{j=0}^{c-d+1} (-1)^j\binom{c-d+1}{j}t^{d+j}. \] 
Substituting this expansion in \eqref{eq:quadratic-Legendre} gives 
\begin{equation*} 
\sum_{k=0}^{n-1} (-1)^k \frac{\binom{n-1}{k}^2}{\binom{n+c}{k+d}}= (n+c+1) \sum_{j=0}^{c-d+1} (-1)^j\binom{c-d+1}{j} \int_0^1t^{d+j}P_{n-1}(1-2t)\,dt. 
\end{equation*} 
By Lemma~\ref{lem:Legendre-moment}, the integral vanishes unless $d+j\geq n-1$. For the remaining terms, Lemma~\ref{lem:Legendre-moment} gives 
\[ \int_0^1 t^{d+j}P_{n-1}(1-2t)\,dt = (-1)^{n-1} \frac{((d+j)!)^2} {(d+j-n+1)!(d+j+n)!}. \] 
This yields \eqref{eq:quadratic-finite-evaluation}. 
\end{proof} 

Finally, taking the special cases of the two sums, we obtain the following interesting identities. \begin{corollary}\label{cor:quadratic-boundary-values} 
Let $c,d$ be integers satisfying $0\leq d\leq c$. Then 
\begin{equation}\label{eq:quadratic-boundary-first} 
\sum_{k=0}^{c+1} (-1)^k \frac{\binom{c+1}{k}^2}{\binom{2c+2}{k+d}} = \frac{(-1)^d}{\binom{2c+2}{c+1}}.
\end{equation} 
Moreover, \begin{equation}\label{eq:quadratic-boundary-second} 
\sum_{k=0}^{c} (-1)^k \frac{\binom{c}{k}^2}{\binom{2c+1}{k+d}} = (-1)^d \frac{c-2d+1}{\binom{2c+1}{c}}.
\end{equation} 
\end{corollary} 
\begin{proof} 
For \eqref{eq:quadratic-boundary-first}, set $n=c+2$ in \eqref{eq:quadratic-finite-evaluation}. The lower and upper limits of the sum on the right are both $c-d+1$, so only one term remains. Straightforward simplification gives 
\[ \sum_{k=0}^{c+1} (-1)^k \frac{\binom{c+1}{k}^2}{\binom{2c+2}{k+d}} = (-1)^d \frac{((c+1)!)^2}{(2c+2)!} = \frac{(-1)^d}{\binom{2c+2}{c+1}}. \] 
For \eqref{eq:quadratic-boundary-second}, set $n=c+1$ in \eqref{eq:quadratic-finite-evaluation}. Only the terms $j=c-d$ and $j=c-d+1$ remain. Their sum simplifies to 
\[ (-1)^d \frac{c-2d+1}{\binom{2c+1}{c}}, \] which proves the result. 
\end{proof} 

Directly from Lemma~\ref{lem:Legendre-moment}, we can derive the following sum.
\begin{proposition}
If $m$ is a non-negative integer and $q$ is a complex number such that $q-m$ is not zero or a negative integer, then
\begin{equation}\label{kx83068}
\sum_{k = 0}^m {( - 1)^k \frac{{\binom{{m}}{k}\binom{{m + k}}{k}}}{{k + q + 1}}}  = \frac{{( - 1)^m q^2 }}{{\left( {q - m} \right)\left( {q + m} \right)\left( {q + m + 1} \right)}}\,\frac{1}{{\left( q \right)_{ - m} \left( q \right)_m }},
\end{equation}
where
\begin{equation*}
\left( x \right)_y  = \frac{{\Gamma \left( {x + y} \right)}}{{\Gamma \left( x \right)}}
\end{equation*}
is the \emph{extended} Pochhammer relation between complex numbers $x$ and $y$ and $\Gamma(x)$ is Euler's Gamma function.
In particular,
\begin{equation}
\sum_{k = 0}^m {( - 1)^k \frac{{\binom{{m}}{k}\binom{{m + k}}{k}}}{{2k + 1}}}  = \frac{1}{{2m + 1}}
\end{equation}
and
\begin{equation}
\sum_{k = 0}^m {( - 1)^k \frac{{\binom{{m}}{k}\binom{{m + k}}{k}}}{{2k + 3}}}  = \frac{1}{{\left( {1 - 2m} \right)\left( {1 + 2m} \right)\left( {3 + 2m} \right)}}.
\end{equation}
\end{proposition}
\begin{proof}
First observe that if $q\ne m$ and $q-m$ is not a negative integer, then~\eqref{eq:Legendre-moment} can be written as
\begin{align}\label{zl2jjb1}
\int_0^1 {t^q \,P_m \left( {1 - 2t} \right)\,dt}  &= \frac{{( - 1)^m q^2 }}{{\left( {q - m} \right)\left( {q + m} \right)\left( {q + m + 1} \right)}}\,\frac{{\Gamma ^2 \left( q \right)}}{{\Gamma \left( {q - m} \right)\Gamma \left( {q + m} \right)}}\nonumber\\
&= \frac{{( - 1)^m q^2 }}{{\left( {q - m} \right)\left( {q + m} \right)\left( {q + m + 1} \right)}}\,\frac{1}{{\left( q \right)_{ - m} \left( q \right)_m }}.
\end{align}
Now plug~\eqref{wc6p34s} into~\eqref{zl2jjb1} after setting $n=m$ and perform the termwise integration to obtain~\eqref{kx83068}.

The particular cases are $q=-\tfrac{1}{2}$ and $q=\tfrac{1}{2}$, respectively.
\end{proof}

\begin{proposition}
If $m$ is a non-negative integer and $q$ is a complex number such that $q-m$ is not zero or a negative integer, then
\begin{align}
&\sum_{k = 0}^m {( - 1)^k \frac{{\binom{{m}}{k}\binom{{m + k}}{k}}}{{\left( {k + q + 1} \right)^2 }}}\nonumber\\
&  = \frac{{( - 1)^m q^2 }}{{\left( {q - m} \right)\left( {q + m} \right)\left( {q + m + 1} \right)}}\,\frac{1}{{\left( q \right)_{ - m} \left( q \right)_m }}\,\left( {H_{q - m - 1}  - 2H_{q - 1}  + H_{q + m - 1} } \right)\nonumber\\
&\qquad\qquad + \frac{(- 1)^m}{(q)_{- m}(q)_m}\,\frac{q\left(2m^3 + 2m^2 + q m^2 + q^3 \right)}{(q - m)^2 (q + m)^2 (q + m + 1)^2}.
\end{align}
\end{proposition}
\begin{proof}
Differentiate~\eqref{kx83068} with respect to $q$ using
\begin{equation*}
\frac{d}{{dx}} \left( x \right)_y = \left( x \right)_y \left( H_{x + y - 1} - H_{x - 1} \right).
\end{equation*}
\end{proof}

\begin{remark}
As we have
$$\binom{n+k}{k}\binom{n}{k}=\binom{n+k}{n}\binom{n}{n-k}=\binom{n+k}{n-k}\binom{2k}{k}=\binom{n+k}{2k}\binom{2k}{k},$$
the above sums can be expressed in a slightly different form. More applications of Legendre polynomials to sums of these forms 
can be found in \cite{Bataille}. For instance, for $0\leq q < n-1$, we have the evaluation \cite[Theorem 3]{Bataille}
\begin{equation}
\sum_{k=0}^{n} \binom{n+k}{2k} \binom{2k}{k} \frac{(-1)^k}{(n+k)(q+k+1)^2} = \frac{(-1)^q (q!)^2 (n-q-2)!}{(n+q+1)!}.
\end{equation}  
\end{remark}

\section{Further applications}

The results stated in the previous sections do not limit to specific sums or polynomial-like identities. For example, by utilizing \eqref{eq:sum_y^jCor23}, they can also be applied to sequences admitting Binet-type representations. As an illustration, we derive an identity valid for every Gibonacci sequence.

So, let $(G_n)_{n\in\mathbb Z}$ be a Gibonacci sequence satisfying \[ G_{n+2}=G_{n+1}+G_n, \] 
with arbitrary initial values $G_0,G_1$. If 
\[ \alpha=\frac{1+\sqrt5}{2}, \qquad \beta=\frac{1-\sqrt5}{2}, \] 
then 
\[ G_n=A\alpha^n+B\beta^n \] 
for suitable constants $A$ and $B$.
\begin{corollary}
The following identity holds
\begin{equation}
  \sum_{k=0}^n\sum_{j=0}^k(-1)^k \frac{\binom{n}{j}}{\binom{n}{k}} G_{m+j} = \frac{(-1)^n((n+1)G_{m+2n}+G_{m+n-2})+(1-n)G_{m+1}+(2n-1)G_m}{n+2}.
\end{equation}
\end{corollary}
\begin{proof}
  It is easy to see that $1+2\alpha=\alpha^3$ and $1+\alpha=\alpha^2$ and similar identities hold for $\beta$. Substitute $y=\alpha$ and $y=\beta$, respectively, to \eqref{eq:sum_y^jCor23} to obtain
  $$\sum_{k=0}^n\sum_{j=0}^k(-1)^k\frac{\binom{n}{j}}{\binom{n}{k}}\alpha^j
		=\frac{1}{n+2}\left((-1)^n(n+1)\alpha^{2n}+\frac{(-1)^n\alpha^{n+2}+n\alpha^2+\alpha^3}{\alpha^4}\right)$$
  and
  $$\sum_{k=0}^n\sum_{j=0}^k(-1)^k\frac{\binom{n}{j}}{\binom{n}{k}}\beta^j
	  =\frac{1}{n+2}\left((-1)^n(n+1)\beta^{2n}+\frac{(-1)^n\beta^{n+2}+n\beta^2+\beta^3}{\beta^4}\right).$$
  Multiplying the first identity by $A\alpha^m$, the second by $B\beta^m$, and adding them, we obtain the result. Note that $G_{m-1}=G_{m+1}-G_m$ and $G_{m-2}=G_m-G_{m-1}=2G_m-G_{m+1}$.    
\end{proof}
A direct specialization to Fibonacci and Lucas numbers gives the following simple identities.
\begin{corollary} 
For every non-negative integer $n$, 
\begin{equation} 
\sum_{k=0}^n\sum_{j=0}^k (-1)^k \frac{\binom{n}{j}}{\binom{n}{k}}F_j = \frac{ (-1)^n\bigl((n+1)F_{2n}+F_{n-2}\bigr)+1-n }{n+2}
\end{equation}
and
\begin{equation} 
\sum_{k=0}^n\sum_{j=0}^k (-1)^k \frac{\binom{n}{j}}{\binom{n}{k}}L_j = \frac{ (-1)^n\bigl((n+1)L_{2n}+L_{n-2}\bigr)+3n-1 }{n+2}. 
\end{equation}
\end{corollary}

We show one more possibillity, which is the following.
\begin{corollary}\label{cor:shifted-weighted-Gibonacci} 
For all non-negative integers $m$ and $n$, 
\begin{align} 
&\sum_{k=0}^n\sum_{j=0}^k (-1)^k \frac{\binom{n}{j}}{\binom{n}{k}} 2^jG_{m+j} \nonumber\\ 
&\quad= \frac{1}{n+2} \left( (-1)^n \bigl((n+1)G_{m+3n} +2^{n+2}G_{m+n-4}\bigr) +nG_{m-3} +4G_{m-5} +G_{m-6} \right). 
\end{align} 
\end{corollary}
\begin{proof}
 Substitute $y=2\alpha$ in \eqref{eq:sum_y^jCor23}, and use 
  $1+2\alpha=\alpha^3$ (and similar expression for $\beta$). We omit the details.
\end{proof}


\begin{thebibliography}{99} \bibitem{Adegoke0}
K. Adegoke, Binomial transforms and the binomial convolution of sequences, arXiv:2507.04179 (2025), 84 pages. \url{https://arxiv.org/abs/2507.04179}

\bibitem{Adegoke3}
K. Adegoke, R. Frontczak and K. Gryszka, Double sums involving binomial coefficients and special numbers, Ann. Math. Sil. (2026), published online October 24, 2026, 33 pages. DOI: 10.2478/amsil-2026-0013.

\bibitem{Adegoke4}
K. Adegoke, R. Frontczak and K. Gryszka, Double sums associated with binomial transforms, arXiv:2511.14799 (2025), 34 pages. \url{https://arxiv.org/abs/2511.14799}

\bibitem{Bataille}
M. Bataille and R. Frontczak, Combinatorial sums derived from properties of Legendre polynomials, Open J. Math. Anal. 10 (1) (2026), 159--169. DOI: 10.30538/psrp-oma2026.0193.

\bibitem{Batir}
N. Batir and A. Sofo, Finite sums involving reciprocals of the binomial and central binomial coefficients and harmonic numbers, Symmetry 13 (11) (2021), Article 2002. DOI: 10.3390/sym13112002.

\bibitem{ChuGuo}
W. Chu and D. Guo, Alternating sums of binomial quotients, Math. Commun. 27 (2) (2022), 203--213.

\bibitem{Gould}
H. W. Gould, {\it Combinatorial Identities}, West Virginia University, Morgantown, 1972.

\bibitem{Graham}
R. L. Graham, D. E. Knuth and O. Patashnik, {\it Concrete Mathematics: A Foundation for Computer Science}, 2nd ed., Addison-Wesley, Reading, MA, 1994.

\bibitem{Gryszka2026}
K. Gryszka, Quadratic and cubic incomplete sums with (non)-symmetric weights, Ramanujan J. 70 (4) (2026), Article 67. DOI: 10.1007/s11139-026-01444-z.

\bibitem{NISTDLMF}
F. W. J. Olver, A. B. Olde Daalhuis, D. W. Lozier, B. I. Schneider, R. F. Boisvert, C. W. Clark, B. R. Miller, B. V. Saunders, H. S. Cohl and M. A. McClain, eds., NIST Digital Library of Mathematical Functions, Release 1.2.8 (2026). \url{https://dlmf.nist.gov/}

\bibitem{OlverEtAl2010}
F. W. J. Olver, D. W. Lozier, R. F. Boisvert and C. W. Clark, eds., {\it NIST Handbook of Mathematical Functions}, Cambridge University Press, New York, 2010.

\bibitem{rathie22}
A. K. Rathie and J. M. Campbell, Generalizations and variants of Knuth's old sum, arXiv:2205.05905 (2022). \url{https://arxiv.org/abs/2205.05905}

\bibitem{Stenlund}
D. Stenlund and J. G. Wan, Some double sums involving ratios of binomial coefficients arising from urn models, J. Integer Sequences 22 (2019), Article 19.1.8.

\bibitem{Szego1975}
G. Szeg\H{o}, {\it Orthogonal Polynomials}, 4th ed., American Mathematical Society Colloquium Publications, Vol. 23, American Mathematical Society, Providence, RI, 1975.

\bibitem{sulei23}
E. Suleiman and B. Sury, A potpourri of pretty identities involving Catalan, Fibonacci and trigonometric numbers, Parabola 59 (2023).

\bibitem{Zhang1}
Z. Zhang, On a kind of curious binomial identity, Discrete Math. 306 (21) (2006), 2740--2754. DOI: 10.1016/j.disc.2006.04.024.

\bibitem{Zhang2}
Z. Zhang and X. Wang, A generalization of Calkin's identity, Discrete Math. 308 (17) (2008), 3992--3997. DOI: 10.1016/j.disc.2007.07.068.
\end{thebibliography}
\end{document}